\documentclass[11pt, a4paper]{amsart}

\usepackage{amsfonts,amsmath,amssymb,amscd}  \usepackage[all]{xy}

\def\N{{\Bbb N}}
\def\Z{{\Bbb Z}}
\def\Q{{\Bbb Q}}

\theoremstyle{definition}

\usepackage{color}

\newcommand\pf{\begin{proof}}
\newcommand\epf{\end{proof}}

\usepackage{mathrsfs}

\numberwithin{equation}{section}

\title[Polynomial values with integer coefficients for the generating functions]
{Polynomial values with integer coefficients for the generating functions \\
of Fibonacci polynomials
}

\author[Yuji~Tsuno]{Yuji~Tsuno}
\address{Yuji Tsuno: National Institute of Technology, Wakayama College,
77 Noshima, Nada-cho, Gobo, Wakayama, Japan 644-0023}
\email{tsuno@wakayama-nct.ac.jp}

\begin{document} 

\begin{abstract}
Fibonacci polynomials are generalizations of Fibonacci numbers, so it is natural to consider polynomial versions of the various results for Fibonacci numbers. 
According to Hong, Pongsriiam, Bulawa, and Lee, the generating function of the Fibonacci sequence in the domain of rational numbers, $f(t)=t/(1-t-t^2)$, takes an integer value if and only if $t=F_{k}/F_{k+1}$ for some $k \in \N$ or $t=-F_{k+1}/F_{k}$ for some $k \in \N^{+}$, where $F_{k}$ is the $k$th Fibonacci number. This study is built upon their work by considering polynomial sequences that satisfy the recurrence relation $F_{i+2}(x)=axF_{i+1}(x)+bF_{i}(x)$ with initial values $(F_{0}(x), F_{1}(x))=(0, 1)$, where $a$ and $b$ are positive integers such that $b|a$. As an application, for a square-free natural number $d \in \N$, we verify the results are of the same form as the above for the generating function of the sequence satisfying the recurrence relation $F_{i+2}(\sqrt{d})=a\sqrt{d} F_{i+1}(\sqrt{d})+bF_{i}(\sqrt{d})$ with initial values $(F_{0}(\sqrt{d}), F_{1}(\sqrt{d}))=(0, 1)$.

\end{abstract}

\maketitle

\noindent

\medskip

\section{Introduction}\label{sec:intro}
Let $a$ and $b$ be positive integers. Bulawa and Lee [1] considered the sequence $\{F_{i}\}_{i \in \N }$ defined by the recurrence relations
\[F_{i+2}=aF_{i+1}+bF_{i}\]
and $F_{0}=0, F_{1}=1$. The generating function is given as
\[f(t)=\frac{t}{1-at-bt^{2}}.\] They established the following necessary and sufficient conditions that should be applied when the rational values in the interval of convergence for the generating function $f(t)$ for the sequence $\{F_{i}\}_{i \in \N }$ are integers:

\vspace{3mm}

\noindent{\bf Theorem 1.1} (Bulawa and Lee [1]).
Let $q$ be a rational number. Let us suppose that $b$ divides $a$ and that $q$ lies within the interval of convergence for the generating function $f(t)$. Then, $f(q) \in \Z$ if and only if 
\[q \in \bigl\{\frac{F_{2i}}{F_{2i+1}}\bigr\}_{i \in \N}.\]

\vspace{3mm}

This result answers the conjecture developed by Hong [2].

\vspace{3mm}

Independent of Bulawa and Lee, Pongsriiam [3] obtained similar results under the conditions of $a = 1$ and $b = 1$.

However, in this study, their work is expanded for applicability to polynomials.

\vspace{3mm}

Let us define a polynomial sequence $\{F_{i}(x)\}_{i \in \N }$ given by the recurrence relations
\[F_{i+2}(x)=axF_{i+1}(x)+bF_{i}(x)\eqno{(1)}\] and $F_{0}(x)=0, F_{1}(x)=1$. For example, $F_{2}(x)=ax$ and $F_{3}(x)=a^{2}x^{2}+b$.
The generating function is given by
\[f(x,t)=\sum_{i=0}^{\infty}F_{i}(x)t^{i}=\frac{t}{1-axt-bt^{2}}.\]
This equation holds within the radius of convergence, but in this paper, we define $f(x, t)=t/(1-axt-bt^{2})$.
Let us also define a polynomial sequence $\{L_{i}(x)\}_{i \in \N }$ given by the recurrence relations
\[L_{i+2}(x)=axL_{i+1}(x)+bL_{i}(x)\eqno{(2)}\] and $L_{0}(x)=2, L_{1}(x)=ax$. For example,  $L_{2}(x)=a^{2}x^{2}+2b, L_{3}(x)=a^{3}x^{3}+3abx$.
The generating function is given by
\[l(x, t)=\sum_{i=0}^{\infty}L_{i}(x)t^{i}=\frac{2-axt}{1-axt-bt^{2}}.\]
This equation holds within the radius of convergence, but in this paper, we define $l(x, t)=(2-axt)/(1-axt-bt^{2})$.
\vspace{2mm}

The main results are as follows:

\vspace{2mm}

\noindent{\bf Theorem 1.2.}
Suppose that $b$ divides $a$
and let $q(x) \in \Q(x)$ be a rational function over $\Q$. For the generating function $f(x, t)$, $f(x, q(x)) \in \Z[x]$ if and only if 
\[q(x) \in \bigl\{\frac{F_{i}(x)}{F_{i+1}(x)}\bigr\}_{i \in \N} \text{ or  }q(x) \in \bigl\{-\frac{F_{i+1}(x)}{bF_{i}(x)}\bigr\}_{i \in \N^{+}}.\]

\vspace{2mm}

\noindent{\bf Theorem 1.3.}
Suppose that $b$ divides $a$ and let $q(x) \in \Q(x)$ be a rational function over $\Q$. For the generating function $l(x, t)$, $l(x, q(x)) \in \Z[x]$ if and only if 
\[q(x) \in \bigl\{\frac{F_{i}(x)}{F_{i+1}(x)}, \frac{L_{i}(x)}{L_{i+1}(x)}, -\frac{L_{i+1}(x)}{bL_{i}(x)}\bigr\}_{i \in \N} \text{ or  }q(x) \in \bigl\{-\frac{F_{i+1}(x)}{bF_{i}(x)}\bigr\}_{i \in \N^{+}}.\]

\vspace{2mm}

Let $d \in \N$ be a square-free natural number. We define a sequence $\{F_{i}(\sqrt{d})\}_{i \in \N }$ given by recurrence relations
\[F_{i+2}(\sqrt{d})=a\sqrt{d} F_{i+1}(\sqrt{d})+bF_{i}(\sqrt{d})\] and $F_{0}(\sqrt{d})=0, F_{1}(\sqrt{d})=1$.
The generating function is given by
\[f(\sqrt{d}, t)=\frac{t}{1-a\sqrt{d} t-bt^{2}}.\]
We define a sequence $\{L_{i}(\sqrt{d})\}_{i \in \N }$ given by recurrence relations
\[L_{i+2}(\sqrt{d})=a\sqrt{d} L_{i+1}(\sqrt{d})+bL_{i}(\sqrt{d})\] and $L_{0}(\sqrt{d})=2, L_{1}(\sqrt{d})=a\sqrt{d}$.
The generating function is given by
\[l(\sqrt{d}, t)=\frac{2-a\sqrt{d}t}{1-a\sqrt{d} t-bt^{2}}.\]

Furthermore, the convergence radii of these generating functions are all 
\[\frac{2}{a\sqrt{d}+\sqrt{a^{2}d+4b}}.\] 
If $d \neq 1$, then in general we do not derive the same result as Theorem 1.2 and Theorem 1.3.

Indeed, we assume that $a=2$, $b=1$ and $d=2$ then
we obtain 
\[ f(\sqrt{d}, \frac{1}{2+\sqrt{d}})=2+\sqrt{d} \in \Z[\sqrt{d}].\]
However, 
\[\frac{1}{2+\sqrt{d}} \notin \bigl\{\frac{F_{i}(\sqrt{d})}{F_{i+1}(\sqrt{d})}\bigr\}_{i \in \N} \text{ and  } \frac{1}{2+\sqrt{d}} \notin \bigl\{-\frac{F_{i+1}(\sqrt{d})}{bF_{i}(\sqrt{d})}\bigr\}_{i \in \N^{+}}.\]
Moreover, $1/(2+\sqrt{d})$ is within the radius of convergence of the generating function $f(\sqrt{d}, t)$.
\vspace{2mm}
In addition, assuming that $a=1, b=1$ and $d=2$ 
we obtain\[ l(\sqrt{d}, \frac{6-5\sqrt{d}}{7})=16-10\sqrt{d} \in \Z[\sqrt{d}].\]
However, 
\[\frac{6-5\sqrt{d}}{7} \notin \bigl\{\frac{F_{i}(\sqrt{d} )}{F_{i+1}(\sqrt{d} )}, \frac{L_{i}(\sqrt{d} )}{L_{i+1}(\sqrt{d} )}, -\frac{L_{i+1}(\sqrt{d} )}{bL_{i}(\sqrt{d} )}\bigr\}_{i \in \N} \text{ and  } \frac{6-5\sqrt{d}}{7}  \notin \bigl\{-\frac{F_{i+1}(\sqrt{d} )}{bF_{i}(\sqrt{d} )}\bigr\}_{i \in \N^{+}}.\]
Moreover, $(6-5\sqrt{d})/7$ is within the radius of convergence of the generating function $l(\sqrt{d}, t)$.
 
If $d=1$, then we have the following theorems.

\noindent{\bf Theorem 1.4.} Suppose that $b$ divides $a$ and $d=1$. Let $q \in \Q$. For the generating function $f(t)$, we have $f(q) \in \Z$ if and only if 
\[q \in \bigl\{\frac{F_{i}}{F_{i+1}}\bigr\}_{i \in \N} \text{ or  } q \in \bigl\{-\frac{F_{i+1}}{bF_{i}}\bigr\}_{i \in \N^{+}},\]where $f(t)=f(\sqrt{d}, t)$ and $F_{i}=F_{i}(\sqrt{d})$. 

\vspace{2mm}

\noindent{\bf Theorem 1.5.}
Suppose that $b$ divides $a$ and $d=1$, and let $q \in \Q$ be a rational number. For the generating function $l(t)$, $l(q) \in \Z$ if and only if 
\[q \in \bigl\{\frac{F_{i}}{F_{i+1}}, \frac{L_{i}}{L_{i+1}}, -\frac{L_{i+1}}{bL_{i}}\bigr\}_{i \in \N} \text{ or  }q \in \bigl\{-\frac{F_{i+1}}{bF_{i}}\bigr\}_{i \in \N^{+}},\]where $f(t)=f(\sqrt{d}, t)$, $F_{i}=F_{i}(\sqrt{d})$ and  $L_{i}=L_{i}(\sqrt{d})$. 

\vspace{2mm}

Focusing on the radii of convergence of the generating functions, we have the following results from Theorem 1.4 and Theorem 1.5.

\vspace{2mm}

\noindent{\bf Corollary 1.6 (= Theorem 1.1)}
Under the assumption of Theorem 1.4, let $q \in \Q$. We assume that $q$ is in the interval of convergence of the generating function $f(t)$. Then, we have $f(q) \in \Z$ if and only if
\[q \in \bigl\{\frac{F_{2i}}{F_{2i+1}}\bigr\}_{i \in \N}.\]

\vspace{2mm}

\noindent{\bf Corollary 1.7.}
Under the assumption of Theorem 1.5, let $q \in \Q$ be a rational number. We assume that $q$ is in the interval of convergence of the generating function $l( t)$. For the generating function $l(t)$, $l(q) \in \Z$ if and only if 
\[q \in \bigl\{\frac{F_{2i}}{F_{2i+1}}, \frac{L_{2i+1}}{L_{2i+2}}\bigr\}_{i \in \N}.\]

\noindent{\bf Remark 1.8.}
Corollary 1.6 and Corollary 1.7 are none other than those given by Bulawa and Lee [1, Theorem 1.1 and Theorem 1.5].

\section{Preliminaries}\label{sec:pre}

Before proceeding to the proof of the main results, the following equations and proposition should be understood.

Let

\[\alpha(x)=\frac{ax+\sqrt{a^{2}x^{2}+4b}}{2},\]
\[
\beta(x)=\frac{ax-\sqrt{a^{2}x^{2}+4b}}{2}
.\]
Then, it follows that
\[F_{n}(x)=\frac{\alpha(x)^{n}-\beta(x)^{n}}{\alpha(x) - \beta(x)}\eqno{(3)}\]
and \[L_{n}(x)=\alpha(x)^{n}+\beta(x)^{n}\eqno{(4)}.\]

Using Equations (3) and (4), the following equations are obtained:

\[F_{n}(x)^{2}-F_{n-1}(x)F_{n+1}(x)=(-b)^{n-1}\eqno{(5)},\]
\[L_{n}(x)^{2}-L_{n-1}(x)L_{n+1}(x)=-(-b)^{n-1}(a^{2}x^{2}+4b)\eqno{(6)},\]
\[F_{2n+1}(x)=L_{n+1}(x)F_{n}(x)+(-b)^{n}\eqno{(7)},\]
\[L_{2n+1}(x)=L_{n+1}(x)L_{n}(x)-(-b)^{n}ax\eqno{(8)},\]
\[L_{2n+1}(x)=(a^{2}x^{2}+4b)F_{n+1}(x)F_{n}(x)+(-b)^{n}ax\eqno{(9)},\]
\[F_{n+1}(x)L_{n}(x)=F_{n}(x)L_{n+1}(x)+2(-b)^{n}\eqno{(10)},\]
\[F_{n+1}(x)=\frac{axF_{n}(x)+L_{n}(x)}{2} \eqno{(11)},\]
\[F_{n}(x)=\frac{-axF_{n+1}(x)+L_{n+1}(x)}{2b} \eqno{(12)},\]
\[L_{n+1}(x)=\frac{axL_{n}(x)+(a^{2}x^{2}+4b)F_{n}(x)}{2} \eqno{(13)},\]
\[L_{n}(x)=\frac{-axL_{n+1}(x)+(a^{2}x^{2}+4b)F_{n+1}(x)}{2b} \eqno{(14)}.\]

The following proposition provides the most robust foundation for the proof of the main results.

\vspace{3mm}

\noindent{\bf Proposition 2.1.}
Let $P(x), Q(x) \in \Q[x]$ be polynomials for which the highest-order coefficient is non-negative. If 
\[
P(x)^{2}-(a^2 x^2 +4b)Q(x)^{2}=4(-b)^{r_{0}} \eqno{(*)}
\]
for some $r_{0} \in \{0,1\}$, then there exists a non-negative integer $n$ such that
\[b^{\lfloor \frac{n}{2}\rfloor}P(x)=L_{n}(x), b^{\lfloor \frac{n}{2}\rfloor}Q(x)=F_{n}(x)\]
and $n \equiv r_{0}$ (mod 2).

\vspace{2mm}

\begin{proof}

First, we define a map of the set of polynomial pairs with rational coefficients satisfying $(*)$ to itself by

\[   \Phi(R(x), S(x))=(\overline{R(x)}, \overline{S(x)})                 \]
where, 
\[   \overline{R(x)}=\frac{-ax(a^{2}x^{2}+4b)S(x)+(a^{2}x^{2}+2b)R(x)}{2b}                          \]
and
\[ \overline{S(x)}=\frac{(a^{2}x^{2}+2b)S(x)-axR(x)}{2b}\]

This is well-defined because $(\overline{R(x)}, \overline{S(x)})$ satisfies the equation $(*)$.
Moreover, the inverse map $\Phi^{-1}$ is given by

\[\Phi^{-1}(R(x), S(x))=(\underline{R(x)}, \underline{S(x)})                 \]
where, 
\[\underline{R(x)}=\frac{ax(a^{2}x^{2}+4b)S(x)+(a^{2}x^{2}+2b)R(x)}{2b}\]
and
\[\underline{S(x)}=\frac{(a^{2}x^{2}+2b)S(x)+axR(x)}{2b}.\]

With the above preparation, to begin, let us consider the case $\mathrm{deg}Q(x) \leq 1.$ 
 
If $Q(x)=0$, we have $P(x)=2$ and $r_{0}=0$. Thus, we obtain 
\[b^{\lfloor \frac{0}{2}\rfloor}P(x)=L_{0}(x) \ \text{and} \ b^{\lfloor \frac{0}{2}\rfloor}Q(x)=F_{0}(x).\]

If $Q(x) \neq 0$ and $\mathrm{deg} Q(x) = 0$, then we have $P(x)=ax$, $Q(x)=1$, and $r_{0}=1$ by matching coefficients of terms with equal degree.
Thus, we obtain
\[b^{\lfloor \frac{1}{2}\rfloor}P(x)=L_{1}(x) \ \text{and} \ b^{\lfloor \frac{1}{2}\rfloor}Q(x)=F_{1}(x).\]

If $\mathrm{deg}Q(x) =1$, the method of undetermined coefficients gives us
\[P(x)=\frac{a^{2}x^{2}+2b}{b}, Q(x)=\frac{ax}{b} \ \ \text{and} \ \ r_{0}=0. \] 
Thus, we obtain the following: \[b^{\lfloor \frac{2}{2}\rfloor}P(x)=L_{2}(x) \ \text{and} \ b^{\lfloor \frac{2}{2}\rfloor}Q(x)=F_{2}(x).\]

Next, let us consider the case $\mathrm{deg}Q(x) \geq 2.$

We have $\mathrm{deg}P(x)=\mathrm{deg}Q(x)+1$ because equation $(*)$  is satisfied.

Let $N=\mathrm{deg}P(x).$ Because equation $(*)$ is satisfied, $c_{0}, c_{1}, \dots, c_{N}$,

$d_{0}, d_{1}, \dots, d_{N-1} \in \Q$ exists
such that $P(x)=c_{0}x^{N}+c_{1}x^{N-1} \dots +c_{N}$

and $Q(x)=d_{0}x^{N-1}+d_{1}x^{N-2} \dots +d_{N-1}$.

Then, we obtain
\[c_{0}=ad_{0}, c_{1}=ad_{1}, ac_{2}=a^{2}d_{2}+2bd_{0}\]
and
\[ac_{3}=\begin{cases}a^{2}d_{3}+2bd_{1} &  \text{if} \ \mathrm{deg}Q(x) > 2 \\ 2bd_{1} &  \text{if} \ \mathrm{deg}Q(x)=2 \end{cases}.\]

On the contrary, if $\mathrm{deg}Q(x) > 2$,
\[\overline{P(x)}=\frac{(a^{2}c_{0}-a^{3}d_{0})x^{N+2}}{2b}+\frac{(a^{2}c_{1}-a^{3}d_{1})x^{N+1}}{2b}+\frac{(a^{2}c_{2}+2bc_{0}-a^{3}d_{2}-4abd_{0})x^{N}}{2b}\]
\[+\frac{(a^{2}c_{3}+2bc_{1}-a^{3}d_{3}-4abd_{1})x^{N-1}}{2b} + \dots.\]
If $\mathrm{deg}Q(x) = 2$, 
\[\overline{P(x)}=\frac{(a^{2}c_{0}-a^{3}d_{0})x^{5}}{2b}+\frac{(a^{2}c_{1}-a^{3}d_{1})x^{4}}{2b}+\frac{(a^{2}c_{2}+2bc_{0}-a^{3}d_{2}-4abd_{0})x^{3}}{2b}\]
\[+\frac{(a^{2}c_{3}+2bc_{1}--4abd_{1})x^{2}}{2b} + (c_{2}-2ad_{2})x+c_{3}.\]
Therefore, we have 
\[\mathrm{deg}\overline{P(x)} \leq \mathrm{deg} P(x)-2\]
from the relationship between the coefficients of $P(x)$ and $Q(x)$.
Moreover, we have
\[\mathrm{deg}\overline{Q(x)} \leq \mathrm{deg} Q(x)-2.\]
Indeed, if $\mathrm{deg}\overline{P(x)}=0$ then $\mathrm{deg}\overline{Q(x)}=0$ since $(\overline{P(x)}, \overline{Q(x)})$ satisfies the equation $(*)$.

If $\mathrm{deg}\overline{P(x)}>0$, we have 
\[\mathrm{deg}\overline{P(x)}=\mathrm{deg}\overline{Q(x)}+1\]
Therefore, \[\mathrm{deg}\overline{Q(x)}=\mathrm{deg}\overline{P(x)}-1 \leq \mathrm{deg}{P(x)}-3=\mathrm{deg}{Q(x)}-2.\]
Alternatively, we have
\[P(x)=\underline{\overline{P(x)}}, Q(x)=\underline{\overline{Q(x)}}\]

Therefore, we see that the highest-order coefficient for $\overline{P(x)}$ (resp. $\overline{Q(x)}$) is a non-negative rational number.
Indeed, if $\mathrm{deg}\overline{P(x)}=0$, $\mathrm{deg}\overline{Q(x)}=0$. Hence, $\overline{P(x)}=\pm 2.$. If $\overline{P(x)}=-2$,   the highest-order coefficient for $P(x)$ is negative. Therefore, $\overline{P(x)}=2$. 
Let $\overline{c_{0}}$ (resp. $\overline{d_{0}}$) be the highest-order coefficient for  $\overline{P(x)}$ (resp. $\overline{Q(x)}$). If $\mathrm{deg}\overline{P(x)} \neq 0$, 
\[\overline{c_{0}}= \pm \overline{d_{0}}a \]
since $(\overline{P(x)}, \overline{Q(x)})$ satisfies the equation $(*)$.
If $\overline{c_{0}}= - \overline{d_{0}}a$, then
\[\mathrm{deg}\overline{P(x)} \nleq \mathrm{deg} P(x)-2.\]
Therefore, we have $\overline{c_{0}} = \overline{d_{0}}a.$ From this, the signs of the highest-order coefficients of $\overline{P(x)}$ and $\overline{Q(x)}$ are equal. If the highest-order coefficient of  $\overline{P(x)}$ is negative, then the highest-order coefficient of  ${P(x)}$ is negative. Therefore, the highest-order coefficient for $\overline{P(x)}$ (resp. $\overline{Q(x)}$) is a non-negative rational number.
Moreover, we have 
\[\mathrm{deg}\overline{P(x)} = \mathrm{deg} P(x)-2\]
and
\[\mathrm{deg}\overline{Q(x)} = \mathrm{deg} Q(x)-2.\]

Writing \[\underbrace{\Phi \circ \dots \circ \Phi}_{n}(P(x), Q(x))=(\overset{n}{\overline{P(x)}}, \overset{n}{\overline{Q(x)}}) \], 
we see that 
\[\mathrm{deg}\overset{ \lfloor \frac{N-1}{2}\rfloor}{\overline{P(x)}} \leq 2.\]

Finally, if there exists a non-negative integer $n$ such that

\[b^{\lfloor \frac{n}{2}\rfloor}\overline{P(x)}=L_{n}(x), b^{\lfloor \frac{n}{2}\rfloor}\overline{Q(x)}=F_{n}(x)\]
and $n  \equiv r_{0}$ (mod 2), then we have

\[b^{\lfloor \frac{n+2}{2}\rfloor}P(x)=L_{n+2}(x), b^{\lfloor \frac{n+2}{2}\rfloor}Q(x)=F_{n+2}(x)\]
and $n+2  \equiv r_{0}$ (mod 2), because
\[L_{m+2}(x)=\frac{ax(a^{2}x^{2}+4b)F_{m}(x)+(a^{2}x^{2}+2b)L_{m}(x)}{2}\]
and
\[F_{m+2}(x)=\frac{(a^{2}x^{2}+2b)F_{m}(x)+axL_{m}(x)}{2}\]
applies to any non-negative integer $m$.
If $\mathrm{deg}P(x) \leq 2$, we have
\[ b^{\lfloor \frac{\mathrm{deg}P(x)}{2}\rfloor}P(x)=L_{\mathrm{deg}P(x)}(x), b^{\lfloor \frac{\mathrm{deg}P(x)}{2}\rfloor}Q(x)=F_{\mathrm{deg}P(x)}(x)                     .\]
Therefore, this completes the proof.

\end{proof}
\vspace{3mm}
The following proposition is obtained by Bulawa and Lee [1, Proposition 1.4], but we prove it by the same proof method as Proposition 2.1.
\vspace{3mm}

\noindent{\bf Proposition 2.2.} Suppose that $b$ divides $a$.
Let $P, Q\in \N$. If 
\[
P^{2}-(a^2 d +4b)Q^{2}=4(-b) \eqno{(**)},
\] then there exists a non-negative integer $n$ such that
\[b^{\lfloor \frac{2n+1}{2}\rfloor}P=L_{2n+1}, b^{\lfloor \frac{2n+1}{2}\rfloor}Q=F_{2n+1},\]where $L_{2n+1}=L_{2n+1}(1), F_{2n+1}=F_{2n+1}(1).$

\vspace{2mm}

\begin{proof}
At first, we define a map of the set of integer pairs satisfying the equation $(**)$ to itself by 
\[\Phi(R, S)=(\overline{R}, \overline{S})                 \]
where, 
\[   \overline{R}=\frac{-a(a^{2}+4b)S+(a^{2}+2b)R}{2b}                          \]
and
\[ \overline{S}=\frac{(a^{2}+2b)S-aR}{2b}\].

This is well-defined as $(\overline{R}, \overline{S})$ satisfies the equation $(**)$.
Moreover, the inverse map $\Phi^{-1}$ is given by

\[   \Phi^{-1}(R, S)=(\underline{R}, \underline{S})                 \]
where, 
\[\underline{R}=\frac{a(a^{2}+4b)S+(a^{2}+2b)R}{2b}\]
and
\[\underline{S}=\frac{(a^{2}+2b)S+aR}{2b}.\]
Under the above preparation,  
to begin, let us consider the case $Q \leq 1$, then we have 
\[P=L_{1}, Q=F_{1}.\]

Next, let us consider the case $Q > 1$, then we have $\overline{P},  \overline{Q} \in \N$. Indeed, we have 
\[P^{2}-a^{2}Q^{2} \in 4\Z\]
since $(P, Q)$ satisfies the equation $(**)$. 
Therefore, $P-aQ \in 2 \Z$. Hence, $\overline{P}, \overline{Q} \in \Z$.
Moreover, $\overline{P}, \overline{Q} \geq 0$.
Indeed, \[\overline{Q}=\frac{(a^{2}+2b)Q-aP}{2b}=\frac{(a^{2}+2b)Q-a\sqrt{(a^{2}+4b)Q^{2}-4b}}{2b}>0\]
and
\[\overline{P}=\frac{-a(a^{2}+4b)Q+(a^{2}+2b)P}{2b}=\frac{-a(a^{2}+4b)Q+(a^{2}+2b)\sqrt{(a^{2}+4b)Q^{2}-4b}}{2b}>0.\]
Moreover, we obtain
\[Q-\overline{Q}=\frac{a\sqrt{(a^{2}+4b)Q^{2}-4b}-a^{2}Q}{2b}>0.\]
In addition, \[Q-\overline{Q} \geq 1,\]
since $Q, \overline{Q} \in \N$.
Writing \[\underbrace{\Phi \circ \dots \circ \Phi}_{n}(P, Q)=(\overset{n}{\overline{P}}, \overset{n}{\overline{Q}}), \] 
we see that there exists a positive integer $m$ such that $\overset{m}{\overline{Q}} \leq 1$
Therefore, the proof is completed in the same way as Proposition 2.1.

\end{proof}

\section{Resulting proofs}

Like [4], the main results are demonstrated by applying Proposition 2.1, Proposition 2.2 and $(1)$ to $(14)$.

\vspace{3mm}
\subsection{Proof of Theorem 1.2.}
First, we suppose that
\[q(x)=\frac{F_{i}(x)}{F_{i+1}(x)} \ (i \in \N)\]
or
\[q(x)=-\frac{F_{i+1}(x)}{bF_{i}(x)} \ (i \in \N^{+}).\]
Then, we show $f(x, q(x)) \in \Z[x]$.
If $i=0$, the result is evident. However, if $i > 0$, by using $(1)$ and $(5)$, we obtain
\[f(x, \frac{F_{i}(x)}{F_{i+1}(x)})=\frac{F_{i}(x)F_{i+1}(x)}{F_{i+1}(x)^{2}-(axF_{i+1}(x)+bF_{i}(x))F_{i}(x)}\]

\[\overset{(1)}=\frac{F_{i}(x)F_{i+1}(x)}{F_{i+1}(x)^{2}-F_{i+2}(x)F_{i}(x)}\overset{(5)}=\frac{F_{i}(x)F_{i+1}(x)}{(-b)^{i}}\] 
\[f(x, -\frac{F_{i+1}(x)}{bF_{i}(x)})=\frac{-bF_{i}(x)F_{i+1}(x)}{bF_{i}(x)(axF_{i+1}(x)+bF_{i}(x))-bF_{i+1}(x)^{2}}\]

\[\overset{(1)}=\frac{-bF_{i}(x)F_{i+1}(x)}{bF_{i}(x)F_{i+2}(x)-bF_{i+1}(x)^{2}}\overset{(5)}=\frac{F_{i}(x)F_{i+1}(x)}{(-b)^{i}}.\]
Moreover, we have $F_{i}(x) \in b^{\lfloor \frac{i}{2}\rfloor}\Z[x] \ (i \in \N)$ 
from the recurrence relation described by Equation (1).  
Indeed, $F_{0}(x) \in b^{\lfloor \frac{0}{2}\rfloor}\Z[x]$ and  $F_{1}(x) \in b^{\lfloor \frac{1}{2}\rfloor}\Z[x]$.
If  $F_{k}(x) \in b^{\lfloor \frac{k}{2}\rfloor}\Z[x]$ and  $F_{k+1}(x) \in b^{\lfloor \frac{k+1}{2}\rfloor}\Z[x]$, then $F_{k+2}(x) \in b^{\lfloor \frac{k+2}{2}\rfloor}\Z[x]$ by Equation (1). Therefore, we have    $F_{i}(x) \in b^{\lfloor \frac{i}{2}\rfloor}\Z[x]$ by mathematical induction. Hence, $f(x, q(x)) \in \Z[x]$.

Next, we suppose that $f(x, q(x))=k(x)$ ($k(x)$ is a polynomial over $\Z$) for some rational function $q(x) \in \Q(x)$, we show that
\[q(x) \in \bigl\{\frac{F_{i}(x)}{F_{i+1}(x)}\bigr\}_{i \in \N} \text{ or  }q(x) \in \bigl\{-\frac{F_{i+1}(x)}{bF_{i}(x)}\bigr\}_{i \in \N^{+}}.\]
If $k(x)=0$, then
\[\frac{q(x)}{1-axq(x)-bq(x)^2}=0.\]
Hence, 
\[q(x)=0=\frac{F_{0}(x)}{F_{1}(x)}\]
If $k(x) \neq 0$, then
\[\frac{q(x)}{1-axq(x)-bq(x)^2}=k(x).\]
Hence,
\[bk(x)q(x)^{2}+(axk(x)+1)q(x)-k(x)=0.\]
Furthermore, 
\[q(x)=\frac{-(axk(x)+1)\pm \sqrt{(axk(x)+1)^2+4bk(x)^2}}{2bk(x)}.\] 
Here, because $q(x)$ is a rational function over $\Q$, there exists a polynomial $M(x) \in \Z[x]$ for which the highest-order coefficient is non-negative such that
\[(axk(x)+1)^{2}+4bk(x)^{2}=M(x)^{2}.\]
This allows us to obtain 
\[\{ (a^{2}x^{2}+4b)k(x)+ax\}^{2}-(a^{2}x^{2}+4b)M(x)^{2}=4(-b).\]
Thus, according to Proposition 2.1, there exists a non-negative integer $n$ such that
\[M(x)=\frac{F_{2n+1}(x)}{b^{n}}, (a^{2}x^{2}+4b)k(x)+ax=\pm \frac{L_{2n+1}(x)}{b^{n}}.\]

From Equation (9), \[(a^{2}x^{2}+4b)k(x)+ax=\pm \frac{(a^{2}x^{2}+4b)F_{n+1}(x)F_{n}(x)+(-b)^{n}ax}{b^n}.\]
Since $k(x) \in \Z[x]$ and $\frac{F_{n}(x)F_{n+1}(x)}{b^{n}} \in \Z[x]$, this means that 
\[(a^{2}x^{2}+4b)k(x)+ax=\frac{L_{2n+1}(x)}{(-b)^{n}}\]
for each $n \in \N^{+}$ given that $k(x) \neq 0$.
Additionally, by using Equation (9) , \[k(x)=\frac{F_{n}(x)F_{n+1}(x)}{(-b)^{n}}.\] 

Consequently, we obtain
\[q(x)=\frac{-axF_{n}(x)F_{n+1}(x)-(-b)^{n}+(-1)^{n}F_{2n+1}(x)}{2bF_{n}(x)F_{n+1}(x)} \ (n \geq 1) \eqno{(A)}\]
or
\[q(x)=\frac{-axF_{n}(x)F_{n+1}(x)-(-b)^{n}-(-1)^{n}F_{2n+1}(x)}{2bF_{n}(x)F_{n+1}(x)} \ (n \geq 1)\eqno{(B)}.\]

Regarding (A) and (B), by using Equations (7), (10), (11), and (12), we obtain
\[q(x) \in \bigl\{\frac{F_{i}(x)}{F_{i+1}(x)}\bigr\}_{i \in \N} \text{ or  }q(x) \in \bigl\{-\frac{F_{i+1}(x)}{bF_{i}(x)}\bigr\}_{i \in \N^{+}}.\]

If $n$ is even, for (A),

\[q(x)=\frac{-axF_{n}(x)F_{n+1}(x)-(-b)^{n}+(-1)^{n}F_{2n+1}(x)}{2bF_{n}(x)F_{n+1}(x)}\overset{(7)}=\frac{-axF_{n+1}(x)+L_{n+1}(x)}{2bF_{n+1}(x)}\]
\[\overset{(12)}=\frac{F_{n}(x)}{F_{n+1}(x)}.\]

If $n$ is odd, for (A),

\[q(x)=\frac{-axF_{n}(x)F_{n+1}(x)-(-b)^{n}+(-1)^{n}F_{2n+1}(x)}{2bF_{n}(x)F_{n+1}(x)}\overset{(7)(10)}=\frac{-axF_{n}(x)-L_{n}(x)}{2bF_{n}(x)}\]
\[\overset{(11)}=-\frac{F_{n+1}(x)}{bF_{n}(x)}.\]

If $n$ is even, for (B),

\[q(x)=\frac{-axF_{n}(x)F_{n+1}(x)-(-b)^{n}-(-1)^{n}F_{2n+1}(x)}{2bF_{n}(x)F_{n+1}(x)}\overset{(7)(10)}=\frac{-axF_{n}(x)-L_{n}(x)}{2bF_{n}(x)}\]
\[\overset{(11)}=-\frac{F_{n+1}(x)}{bF_{n}(x)}.\]

If $n$ is odd, for (B),

\[q(x)=\frac{-axF_{n}(x)F_{n+1}(x)-(-b)^{n}-(-1)^{n}F_{2n+1}(x)}{2bF_{n}(x)F_{n+1}(x)}\overset{(7)}=\frac{-axF_{n+1}(x)+L_{n+1}(x)}{2bF_{n+1}(x)}\]
\[\overset{(12)}=\frac{F_{n}(x)}{F_{n+1}(x)}.\]

\vspace{3mm}
\subsection{Proof of Theorem 1.3.} First, we suppose that 
 \[q(x)=\frac{F_{i}(x)}{F_{i+1}(x)}, \frac{L_{i}(x)}{L_{i+1}(x)}, -\frac{L_{i+1}(x)}{bL_{i}(x)} \ (i \in \N) \text{ or  }  q(x)=-\frac{F_{i+1}(x)}{bF_{i}(x)} \  (i \in \N^{+}),\]then, we show that $l(x, q(x)) \in \Z[x]$

If $i=0$, the result is evident. However, if $i > 0$, by using Equations $(1), (2), (5), (6), (13)$ and $(14)$, we obtain
\[l(x, \frac{F_{i}(x)}{F_{i+1}(x)})=\frac{2F_{i+1}(x)^{2}-axF_{i}(x)F_{i+1}(x)}{F_{i+1}(x)^{2}-F_{i}(x)(axF_{i+1}(x)+bF_{i}(x))}\overset{(1)}=\frac{2F_{i+1}(x)^{2}-axF_{i}(x)F_{i+1}(x)}{F_{i+1}(x)^{2}-F_{i}(x)F_{i+2}(x)}\]
\[\overset{(5)}=\frac{2F_{i+1}(x)^{2}-axF_{i}(x)F_{i+1}(x)}{(-b)^{i}}\]
\[l(x, \frac{L_{i}(x)}{L_{i+1}(x)})=\frac{2L_{i+1}(x)^{2}-axL_{i}(x)L_{i+1}(x)}{L_{i+1}(x)^{2}-L_{i}(x)(axL_{i+1}(x)+bL_{i}(x))}\overset{(2)}=\frac{2L_{i+1}(x)^{2}-axL_{i}(x)L_{i+1}(x)}{L_{i+1}(x)^{2}-L_{i}(x)L_{i+2}(x)}\]
\[\overset{(6)}=\frac{L_{i+1}(x)(2L_{i+1}(x)-axL_{i}(x))}{-(-b)^{i}(a^{2}x^{2}+4b)}\overset{(13)}=\frac{L_{i+1}(x)F_{i}(x)}{-(-b)^{i}} \]
\[l(x, -\frac{L_{i+1}(x)}{bL_{i}(x)})=\frac{L_{i}(x)(2bL_{i}(x)+axL_{i+1}(x))}{L_{i}(x)(axL_{i+1}(x)+bL_{i}(x))-L_{i+1}(x)^{2}} \ \ \ \ \ \ \ \ \ \ \ \ \ \ \ \ \ \ \ \ \ \ \ \ \ \ \  \]
\[\overset{(2)}=\frac{L_{i}(x)(2bL_{i}(x)+axL_{i+1}(x))}{L_{i}(x)L_{i+2}(x)-L_{i+1}(x)^{2}}\overset{(6)(14)}=\frac{L_{i}(x)F_{i+1}(x)}{(-b)^{i}}\]
\[l(x, -\frac{F_{i+1}(x)}{bF_{i}(x)})=\frac{F_{i}(x)(2bF_{i}(x)+axF_{i+1}(x))}{F_{i}(x)(bF_{i}(x)+axF_{i+1}(x))-F_{i+1}(x)^{2}}\ \ \ \ \ \ \ \ \ \ \ \ \ \ \ \ \ \ \ \ \ \ \ \ \ \ \  \]
\[\overset{(1)(5)}=\frac{F_{i}(x)(2bF_{i}(x)+axF_{i+1}(x))}{F_{i}(x)F_{i+2}(x)-F_{i+1}(x)^{2}}=\frac{F_{i}(x)(2bF_{i}(x)+axF_{i+1}(x))}{-(-b)^{i}}.\]

Using an induction similar to the discussion in the proof of Theorem 1.2, we obtain that $F_{i}(x), L_{i}(x) \in b^{\lfloor \frac{i}{2}\rfloor}\Z[x]$ from the recurrence relations described by Equations (1) and (2),
\[l(x, \frac{F_{i}(x)}{F_{i+1}(x)}), l(x, \frac{L_{i}(x)}{L_{i+1}(x)}), l(x, -\frac{L_{i+1}(x)}{bL_{i}(x)}) \in \Z[x] \ (i \in \N) \text{ and  }  l(x, -\frac{F_{i+1}(x)}{bF_{i}(x)}) \in \Z[x] \  (i \in \N^{+}).\]

Next, we suppose that $l(x, q(x))=k(x)$ ($k(x)$ is a polynomial over $\Z$) for some rational function $q(x) \in \Q(x)$, then we show that
\[q(x) \in \bigl\{\frac{F_{i}(x)}{F_{i+1}(x)}, \frac{L_{i}(x)}{L_{i+1}(x)}, -\frac{L_{i+1}(x)}{bL_{i}(x)}\bigr\}_{i \in \N} \text{ or  }q(x) \in \bigl\{-\frac{F_{i+1}(x)}{bF_{i}(x)}\bigr\}_{i \in \N^{+}}.\]

If $k(x)=0$, then

\[\frac{2-axq(x)}{1-axq(x)-bq(x)^2}=0.\]

Therefore, 
\[q(x)=\frac{2}{ax}=\frac{L_{0}(x)}{L_{1}(x)}.\]
Alternatively, if $k(x) \neq 0$, then
\[\frac{2-axq(x)}{1-axq(x)-bq(x)^2}=k(x).\]
Hence,
\[bk(x)q(x)^{2}+ax(k(x)-1)q(x)+2-k(x)=0.\]
Therefore, 
\[q(x)=\frac{-ax(k(x)-1) \pm \sqrt{a^{2}x^{2}(k(x)-1)^{2}-4bk(x)(2-k(x))}}{2bk(x)}.\] 
Here, because $q(x)$ is a rational function over $\Q$, there exists a polynomial $M(x) \in \Z[x]$ for which the highest-order coefficient is non-negative such that
\[a^{2}x^{2}(k(x)-1)^{2}-4bk(x)(2-k(x))=M(x)^{2}.\]
Then, we obtain 
\[M(x)^{2}-(a^{2}x^{2}+4b)(k(x)-1)^{2}=4(-b).\]
Thus, according to Proposition 2.1, there exists a non-negative integer $n$ such that
\[M(x)=\frac{L_{2n+1}(x)}{b^{n}}, k(x)-1=\pm \frac{F_{2n+1}(x)}{b^{n}}.\]

Hence, we obtain
\[q(x)=\frac{-axF_{2n+1}(x)+L_{2n+1}(x)}{2b(F_{2n+1}(x)+b^{n})} \ (n \geq 0)\eqno{(C)},\]

\[q(x)=\frac{-axF_{2n+1}(x)-L_{2n+1}(x)}{2b(F_{2n+1}(x)+b^{n})} \ (n \geq 0)\eqno{(D)},\]

\[q(x)=\frac{axF_{2n+1}(x)+L_{2n+1}(x)}{2b(-F_{2n+1}(x)+b^{n})} \ (n \geq 1) \eqno{(E)},\]
or
\[q(x)=\frac{axF_{2n+1}(x)-L_{2n+1}(x)}{2b(-F_{2n+1}(x)+b^{n})} \ (n \geq 1) \eqno{(F)}.\]

In the cases of (C) to (F), using Equations (7), (9), (10), (13), and (14) gives us
\[q(x) \in \bigl\{\frac{F_{i}(x)}{F_{i+1}(x)}, \frac{L_{i}(x)}{L_{i+1}(x)}, -\frac{L_{i+1}(x)}{bL_{i}(x)}\bigr\}_{i \in \N} \text{ or  }q(x) \in \bigl\{-\frac{F_{i+1}(x)}{bF_{i}(x)}\bigr\}_{i \in \N^{+}}.\]

Then, if $n$ is even, for (C), 

\[q(x)=\frac{-axF_{2n+1}(x)+L_{2n+1}(x)}{2b(F_{2n+1}(x)+b^{n})}\overset{(7)(9)(10)}=\frac{-axL_{n+1}(x)F_{n}(x)+(a^{2}x^{2}+4b)F_{n+1}(x)F_{n}(x)}{2bF_{n+1}(x)L_{n}(x)}\]
\[\overset{(14)}=\frac{F_{n}(x)}{F_{n+1}(x)}.\]

If $n$ is odd, for (C), 

\[q(x)=\frac{-axF_{2n+1}(x)+L_{2n+1}(x)}{2b(F_{2n+1}(x)+b^{n})}\overset{(7)(9)}=\frac{-axL_{n+1}(x)F_{n}(x)+(a^{2}x^{2}+4b)F_{n+1}(x)F_{n}(x)}{2bL_{n+1}(x)F_{n}(x)}\]
\[\overset{(14)}=\frac{L_{n}(x)}{L_{n+1}(x)}.\]

If $n$ is even, for (D),

\[q(x)=\frac{-axF_{2n+1}(x)-L_{2n+1}(x)}{2b(F_{2n+1}(x)+b^{n})}\overset{(7)(9)(10)}=\frac{-axL_{n}(x)-(a^{2}x^{2}+4b)F_{n}(x)}{2bL_{n}(x)}\]

\[\overset{(13)}=-\frac{L_{n+1}(x)}{bL_{n}(x)}.\]

If $n$ is odd, for (D), 

\[q(x)=\frac{-axF_{2n+1}(x)-L_{2n+1}(x)}{2b(F_{2n+1}(x)+b^{n})}\overset{(7)(9)(10)}=\frac{-axF_{n+1}(x)L_{n}(x)-(a^{2}x^{2}+4b)F_{n+1}(x)F_{n}(x)}{2bL_{n+1}(x)F_{n}(x)}\]

\[\overset{(13)}=-\frac{F_{n+1}(x)}{bF_{n}(x)}.\]

If $n$ is even, for (E),

\[q(x)=\frac{axF_{2n+1}(x)+L_{2n+1}(x)}{2b(-F_{2n+1}(x)+b^{n})}\overset{(7)(9)(10)}=\frac{axF_{n+1}(x)L_{n}(x)+(a^{2}x^{2}+4b)F_{n+1}(x)F_{n}(x)}{-2bL_{n+1}(x)F_{n}(x)}\]

\[\overset{(13)}=-\frac{F_{n+1}(x)}{bF_{n}(x)}.\]

If $n$ is odd, for (E), 

\[q(x)=\frac{axF_{2n+1}(x)+L_{2n+1}(x)}{2b(-F_{2n+1}(x)+b^{n})}\overset{(7)(9)(10)}=\frac{axF_{n+1}(x)L_{n}(x)+(a^{2}x^{2}+4b)F_{n+1}(x)F_{n}(x)}{-2bF_{n+1}(x)L_{n}(x)}\]

\[\overset{(13)}=-\frac{L_{n+1}(x)}{bL_{n}(x)}.\]

If $n$ is even, for (F),

\[q(x)=\frac{axF_{2n+1}(x)-L_{2n+1}(x)}{2b(-F_{2n+1}(x)+b^{n})}\overset{(7)(9)}=\frac{axL_{n+1}(x)F_{n}(x)-(a^{2}x^{2}+4b)F_{n+1}(x)F_{n}(x)}{-2bL_{n+1}(x)F_{n}(x)}\]

\[\overset{(14)}=\frac{L_{n}(x)}{L_{n+1}(x)}.\]

If $n$ is odd, for (F), 

\[q(x)=\frac{axF_{2n+1}(x)-L_{2n+1}(x)}{2b(-F_{2n+1}(x)+b^{n})}\overset{(7)(9)(10)}=\frac{axL_{n+1}(x)F_{n}(x)-(a^{2}x^{2}+4b)F_{n+1}(x)F_{n}(x)}{-2bF_{n+1}(x)L_{n}(x)}\]

\[\overset{(14)}=\frac{F_{n}(x)}{F_{n+1}(x )}.\]

\vspace{3mm}
\subsection{Proof of Theorem 1.4.}

In the proof of Theorem 1.2, we apply Proposition 2.2 instead of applying Proposition 2.1 and set $x = 1$ to complete the proof.

\vspace{3mm}
\subsection{Proof of Theorem 1.5.}
In the proof of Theorem 1.3, we apply Proposition 2.2 instead of applying Proposition 2.1 and set $x = 1$ to complete the proof.

\section{Acknowledgement}
I thank Editage for the English language editing.

\vspace{3mm}

\noindent
{\sc Mathematics Subject Classification (2010):}
11B39 . 


\begin{thebibliography}{99}

\bibitem{BL}
A.~Bulawa and W.~K.~Lee,\
\emph{Integral values of the generating functions for the Fibonacci and related sequences},
The Fibonacci Quarterly, \textbf{55.1} (2017), 74--81.
\bibitem{Ho}
D.~S.~Hong,\ 
\emph{When is the generating function of the Fibonacci numbers an integer?},
The College Mathematics Journal, \textbf{46} (2015), 110--112.

\bibitem{Po}
P.~Pongsriiam,\
\emph{Integral values of the generating functions of Fibonacci and Lucas numbers},
The College Mathematics Journal, \textbf{48} (2017), 97--101.

\bibitem{Tsu}
Y.~Tsuno,\
\emph{Extended results on integer values of generating functions for sequences given by Pell's equations},
The Fibonacci Quarterly, \textbf{59.2} (2021), 158--166.

\end{thebibliography}
\end{document}